\documentclass[10pt,twocolumn]{article}

\usepackage{graphicx}%
\usepackage{color}%
\usepackage{amsmath,amssymb,amsthm,eucal,upref,bm}%
\usepackage{mathtools}%
\usepackage{txfonts}%
\usepackage{mathptmx}%

\usepackage{fancybox}

\usepackage{mathrsfs}%

\usepackage{wrapfig} 

\usepackage{sidecap}

\usepackage{setspace}

\usepackage{ulem}

\usepackage[small,compact]{titlesec}%
\usepackage[small,it]{caption}

\allowdisplaybreaks[1]

\makeatletter

\usepackage[paperwidth=210mm,paperheight=297mm,top=22mm,bottom=22mm,left=20mm,right=20mm]{geometry}

\let\leq=\leqslant\let\geq=\geqslant

\def\({{\rm(}} \def\){\/{\rm)}} 

\newcommand{\D}{\mathrm{d}\kern0.2pt}% 
\newcommand{\RR}{\mathbb{R}}% 
\newcommand{\CC}{\mathbb{C}}%

\def\comm#1{\mbox{{\large\ding{45}}}%
\setbox0=\hbox to \textwidth{\kern-16mm\raise0.6ex\hbox{\ovalbox{%
\vtop{\tolerance=7000\hyphenpenalty=3000\leftskip=0pt\rightskip=0pt plus 1pt\relax
\hsize=13mm\tiny\noindent #1}}}\hfill}%
\dp0=0pt\ht0=0pt\vadjust{\box0}}%

\begin{document}

\renewcommand{\topfraction}{0.9} 
\renewcommand{\bottomfraction}{0.9}
\renewcommand{\thefootnote}{}

\abovedisplayshortskip=-2pt plus 0.25pt \belowdisplayshortskip=7pt plus 0.25pt minus 1pt
\abovedisplayskip=9.5pt plus 0.25pt minus 1pt \belowdisplayskip=9.5pt plus 0.25pt minus
1pt

%\baselineskip=0.0188\textheight

%\lineskip=0cm plus1pt \hfuzz=4pt

\linespread{1}

\title{{\bfseries {\fontsize{14}{\f@baselineskip}\selectfont The legacy of J\'ozef
Marcinkiewicz: four hallmarks of genius}}}

\author{Nikolay Kuznetsov}

\date{In memoriam of an extraordinary analyst}

\twocolumn[\maketitle]

\setcounter{equation}{0}

This article is a tribute to one of the most prominent Polish mathematicians
J\'o\-zef Marcinkiewicz who perished 80 years ago in the Katy\'n massacre. He was
one of nearly 22000 Polish officers interned by the Red Army in September 1939 and
executed in April--May 1940 in the Katy\'n forest near Smolensk and at several
locations elsewhere. One of these places was Kharkov (Ukraine), where more than 3800
Polish prisoners of war from the Starobelsk camp were executed. One of them was
Marcinkiewicz; the plaque with his name (see photo at the bottom) is on the Memorial
Wall at the Polish War Cemetery in Kharkov.$^*$ \footnote{$^*$
https://www.tracesofwar.com/sights/10355/Polish-War-Cemetery-Kharkiv.htm} This
industrial execution was authorized by Stalin's secret order dated 5 March 1940,
and organised by Beria, who headed the People's Commissariat for Internal Affairs
(the interior ministry of the Soviet Union) known as NKVD.

Turning to the personality and mathematical achievements of Marcinkiewicz, it is
appropriate to cite the article \cite{Z} of his superviser Antoni Zygmund (it is
published in the {\it Collected Papers} \cite{JM} of Marcinkiewicz;
see~p.~1):\\[-5mm]
\begin{quote}
Considering what he did during his short life and what he might have done in normal
circumstances one may view his early death as a great blow to Polish Mathematics,
and probably its heaviest individual loss during the second world war.\\[-4mm]
\end{quote}

\begin{figure}[b]
\centering
\vspace{-2mm}
\includegraphics[width=64mm]{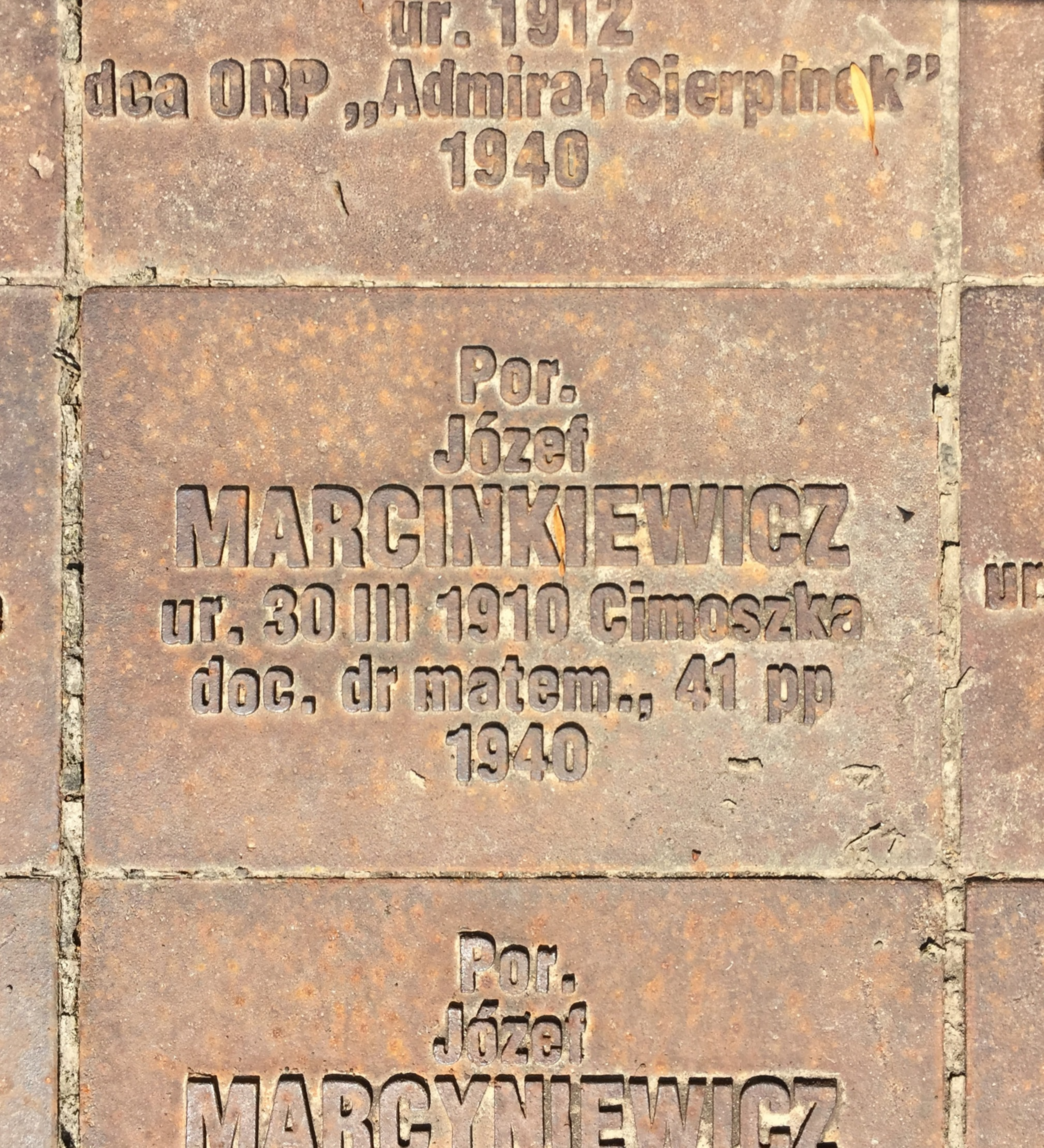}
\end{figure}

\subsection*{From the Marcinkiewicz biography \cite{LM1}}

On the occasion of the centenary of Marcinkiewicz's birth, a conference was held on
28 June--2 July 2010, in Pozna\'n. In its proceedings, L.~Maligranda published the
detailed article \cite{LM1} about Marcinkiewicz's life and mathematical results; 16
pages of this paper are devoted to his biography, where one finds the following
about his education and scientific career.

{\bf Education.} Klemens Marcinkiewicz, J\'ozef's father, was a well-to-do farmer to
afford private lessons for him at home (the reason was J\'ozef's poor health),
before sending him to elementary school and then to gymnasium in Bia\l y\-stok.
After graduating from it in 1930, J\'ozef enrolled in the Department of Mathematics
and Natural Science of the Stefan Batory University (USB) in Wilno (then in Poland,
now Vilnius in Lithuania).

\begin{figure}[t]
\centering
\includegraphics[width=64mm]{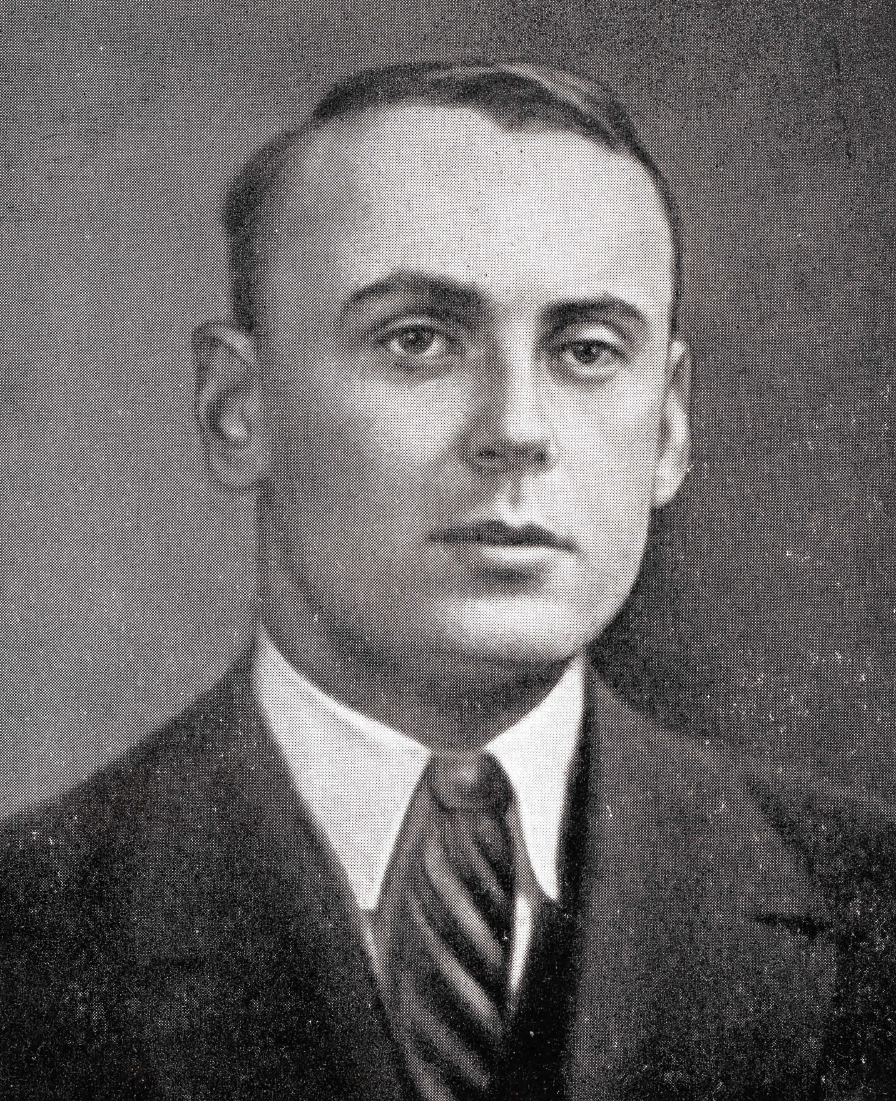}
\vspace{-5mm}
\end{figure}

From the beginning of his university studies, J\'ozef de\-monstrated exceptional
mathematical talent which attracted attention of his professors, in particular, of
A. Zygmund. His lectures on orthogonal series, requiring some erudition, in
particular, knowledge of the Lebesgue integral, Marcinkiewicz attended being just a
second year student; this was the point, where their collaboration began. The first
paper of Marcinkiewicz (see \cite{JM}, p.~35) had been published when he was still
an undergraduate student. It provides a half-page proof of Kolmogorov's theorem
(1924) guaranteeing the convergence almost everywhere for partial sums of lacunary
Fourier series. Marcinkiewicz completed his MSc and PhD theses (both supervised by
Zygmund) in 1933 and 1935 respectively; to obtain PhD degree he also passed a rather
stiff examination. The second dissertation was the fourth of almost five dozens his
publications; it concerns interpolation by means of trigonometric polynomials and
contains interesting results (see \cite{Z}, p.~17, for a discussion), but a long
publication history awaited this work. Its part was published in the {\it Studia
Mathematica} the next year after the thesis defence (these two papers in French are
reproduced in \cite{JM}, pp.~171--185 and 186--199). The full, original text in
Polish appeared in the {\it Wiadomo\'sci Matematyczne} (the {\it Mathematical News})
in 1939; finally, its English translation was included in \cite{JM}, pp.~45--70.

{\bf Scientific career.} During the two years between defending his MSc and PhD
theses, Marcinkiewicz did the one year of mandatory military service and then was
Zygmund's assistant at USB. The aca\-demic year 1935/1936, Marcinkiewicz spent as an
assis\-tant at the Jan Kazimierz University in Lw\'ow. Despite 12 hours of teaching
weekly, he was an active participant of mathematical discussions at the famous
Scottish Caf\'e (see \cite{D}, ch.\,10, where this unique form of doing mathematics
is described), and his contribution to the {\it Scottish Book} compiled in this
caf\'e was substantial taking into account that his stay in Lw\'ow lasted only nine
months. One finds the history of this book in \cite{SB}, ch.\,I, whereas problems
and their solutions, where applicable, are presented in ch.\,II. Marcinkiewicz posed
his own problem; it concerns the uniqueness of solution for the integral equation
\[ \int_0^1 y (t) f (x-t) \, \D t = 0 , \ \ x \in [0, 1] ;
\]
he conjectured that if $f (0) \neq 0$ and $f$ is continuous, then this equation has
only the trivial solution $y \equiv 0$ (see problem no. 124, \cite{SB}, pp.~211 and
212). He also solved three problems; his negative answers to problems 83 and 106
posed by H. Auerbach and S.~Banach, respectively, involve ingenious counterexamples.
His positive solution of problem 131 (it was formulated by Zygmund in a lecture
given in Lw\'ow in the early 1930s) was published in 1938; see \cite{JM},
pp.\,413--417.

During the next two academic years, Marcinkiewicz was a senior assistant at USB and
after completing his habilitation in June 1937 became the youngest docent at USB.
The same year, he was awarded the J\'ozef Pi\l sudski Scientific Prize (the highest
Polish distinction for achievements in science at that time). His last academic year
1938/1939, Marcinkiewicz was on leave from USB; a scholarship from the Polish Fund
for National Culture yielded him opportunity to travel. He spent October 1938--March
1939 in Paris and moved to the University College London for April--Au\-gust 1939,
visiting also Cambridge and Oxford.

This period was very successful for Marcinkiewicz; he published several brief notes
in the {\it Comptes rendus de l'Acad\'emie des Sciences Paris}. One of these,
namely~\cite{M3}, became widely cited because the celebrated theorem concerning
interpolation of operators was announced in it. Now, this theorem is referred to as
the Marcinkiewicz or Marcinkiewicz--Zygmund interpolation theorem (see below).
Moreover, an important notion was introduced in the same note; the so-called
weak-$L^p$ spaces, known as Marcinkiewicz spaces now, are essential for the general
form of this theorem.

Meanwhile, Marcinkiewicz was appointed to the position of Extraordinary Professor at
the University of Pozna\'n in June 1939. On his way to Paris, he delivered a lecture
there and this, probably, was related to this impending appointment; also, this was
the reason to decline an offer of professorship in the USA during his stay in Paris.

Marcinkiewicz still was in England, when the general mobilisation was announced in
Poland in the second half of August 1939; the outbreak of war became imminent. His
colleagues advised him to stay in England, but his ill-fated decision was to go back
to Poland. He regarded himself a patriot of his homeland, which is easily
explainable by the fact that he was just eight years old (very sensitive age in
forming a personality) when the independence of Poland was restored.

\subsection*{Contribution of Marcinkiewicz to mathematics}

Marcinkiewicz was a prolific author as demonstrates a list of his almost five dozen
papers written just in seven years (1933--1939); see {\it Collected Papers}
\cite{JM}, pp.~31--33. He was open to collaboration; indeed, more than one third of
his papers (19, to be exact) were written with five coauthors, of which the lion's
share belongs to his superviser Zygmund.

Marcinkiewicz is known, primarily, as an outstanding analyst, whose best results
deal with various aspects of real analysis; in particular, theory of series
(trigonometric and others), inequalities and approximation theory. He also published
several papers concerning complex and functional analysis and probability theory.
In the extensive paper \cite{LM1} dedicated to the centenary of Marcinkiewicz's
birth, one finds a detailed survey of all his results.

This survey begins with the description of five topics concerning functional
analysis (\cite{LM1}, pp.~153--175). No doubt, the first two of them---the
Marcinkiewicz interpolation theorem and Marcinkiewicz spaces---are hallmarks of
genius. An indirect evidence of ingenuity of the idea behind these results is that
the note \cite{M1}, in which they first appeared, is the most cited work of
Marcinkiewicz.

Another important point about his work is that he skillfully applied methods of real
analysis to questions bordering with complex analysis. A brilliant example of this
mas\-tery---one more hallmark of genius---is the Marcinkiewicz function $\mu$
introduced as an analogue of the Littlewood--Paley function $g$. It is worth
mentioning that the short paper \cite{M4}, in which $\mu$ first appeared, contains
other fruitful ideas developed by many mathematicians subsequently.

One more hallmark of genius one finds in the paper \cite{M1} entitled ``Sur les
multiplicateurs des s\'eries de Fourier''. There are many generalisations of its
results because of their important applications. This work was the last of eight
papers that Marcinkiewicz published in the {\it Studia Mathematica}; the first three
he submitted during his stay in Lw\'ow and they appeared in 1936.

Below, the above mentioned results of Marcinkiewicz are outlined in their historical
context together with some further developments. One can find a detailed
presentation of all these results in the excellent textbook \cite{S} based on
lectures of the eminent analyst Elias Stein, who made a considerable contribution to
further development of ideas proposed by Marcinkiewicz.

\subsection*{Marcinkiewicz interpolation theorem \\ and Marcinkiewicz spaces}

There are two pillars of the interpolation theory: the classical Riesz--Thorin and
Marcinkiewicz theorems. Each of these serves as the basis for two essentially
different approaches to interpolation of operators known as the complex and real
methods. The term `interpolation of operators' was, presumably, coined by
Marcinkiewicz in 1939, because Riesz and Thorin, who published their results in 1926
and 1938, respectively, referred to their assertions as `convexity theorems'.

It is worth emphasising again that a characteristic feature of Mar\-cin\-kiewicz's
work was applying real methods to problems that other authors treated with the help
of complex analysis. It was mentioned above that in his paper~\cite{M4} published in
1938, Marcinkiewicz introduced the function $\mu$ without using complex variables,
but so that it is analogous to the Littlewood--Paley function $g$, whose definition
involves these variables. In the same year 1938, Thorin published his extension of
the Riesz convexity theorem, which exemplifies the approach based on complex
variables. Possibly, this stimulated Marcinkiewicz to seek an analogous result with
proof relying on real analysis. Anyway, Marcinkiewicz found his interpolation
theorem and announced in \cite{M3}; concurrently, a letter was sent to Zygmund which
contained the proof concerning a particular case. Ten years after World War II,
Zygmund reconstructed the general proof and published it in 1956; for this reason
the theorem is sometimes referred to as the Marcinkiewicz--Zygmund interpolation
theorem.

An excellent introduction to the interpolation theory one finds in the book
\cite{BL} based on the works of Jaak Peetre (he passed away on 1 April 2019 aged
83), whose contribution to this theory cannot be overestimated. In collaboration
with Jacques-Louis Lions, he introduced the `real method interpolation spaces' (see
their fundamental article \cite{LP}), which can be considered as `descendants' of
the Marcinkiewicz interpolation theorem.

An important fact of Peetre's biography is that his life was severely changed during
World War II (another reminder about that terrible time). With his parents, Jaak
escaped from Estonia in September 1944 just two days before his home town P\"arnu
was destroyed in an air raid of the Red Army. He was only ten years old when his
family settled in Lund (Sweden), where he spent most of his life. But let us turn to
mathematics again.

{\bf The Marcinkiewicz interpolation theorem for operators in $\mathbf{L^p
(\RR^n)}$.} We begin with this simple result because it has numerous applications,
being valid for sub-additive operators mapping the Lebesgue spaces $L^p (\RR^n)$
with $p \geq 1$ into themselves (see, e.g., \cite{S}, ch.~1, sect.~4). We recall
that an operator $T: L^p \to L^p$ is sub-additive if
\[ | T(f_1 + f_2) (x)| = | T(f_1) (x) | + | T(f_2) (x) | \ \ \mbox{for every} \ f_1, 
f_2 \, .
\]
Furthermore, $T$ is of {\it weak type} $(r,r)$ if the inequality
\[ \alpha^r {\rm mes} \{ x : | T(f) (x) | > \alpha \} \leq A_r \| f \|_r^r
\]
holds for all $\alpha > 0$ and all $f \in L^r$ with $A_r$ independent of $\alpha$
and $f$; here, ${\rm mes} \{ \dots \}$ denotes the Lebesgue measure of the
corresponding set and
\[ \| f \|_p = \Big[ \int_{\RR^n} | f (x) |^p \, \D x \Big]^{1/p}
\]
is the norm in $L^p (\RR^n)$. Now, we are in a position to formulate the following.

\vspace{1mm}

\noindent {\bf Theorem 1.} {\it Let $1 \leq r_1 < r_2 < \infty$, and let $T$ be a
sub-additive operator acting simultaneously in $L^{r_i} (\RR^n)$, $i = 1, 2$. If it
is of weak type $(r_i, r_i)$ for $i = 1, 2$, then for every $p \in (r_1, r_2)$ the
inequality $\| T (f) \|_p \leq B \| f \|_p$  holds for all $f \in L^p (\RR^n)$ with
$B$ depending only on $A_{r_1},\ A_{r_2},\ r_1,\ r_2$ and $p$.}

\vspace{1mm}

When $B$ is independent of $f$ in the last inequality, the operator $T$ is of {\it
strong type} $(p, p)$; it is clear that $T$ is also of weak type $(p, p)$ in this
case.

In the letter to Zygmund mentioned above, Marcin\-kie\-wicz included a proof of
this theorem for the case $r_1 = 1$ and $r_2 = 2$. Presumably, it was rather simple;
indeed, even when $r_2 < \infty$ is arbitrary, the proof is less than two pages long
in \cite{S}, ch.~1, sect.~4.

{\bf Marcinkiewicz spaces.} Another crucial step, made by Marcinkiewicz in
\cite{M3}, was introduction of the {\it weak} $L^p$ spaces playing the essential
role in his general interpolation theorem. Now, they are called the {\it
Mar\-cin\-kiewicz spaces} and usually denoted $L^{p,\infty}$.

To give an idea of these spaces, let us consider a measure space $(U, \Sigma, m)$
over real scalars with a non-negative measure $m$ (just to be specific). For a
real-valued $f$, which is finite almost everywhere and $m$-measurable, we introduce
its distribution function\\[-2mm]
\[ m ( \{ x : | f (x) | > \lambda \} ) , \ \ \lambda \in (0, \infty) ,\\[-1mm]
\]
and put\\[-2mm]
\[ | f |_{p,\infty} = \sup_{\lambda > 0} \, \lambda [ m ( \{ x : | f (x) | > \lambda 
\} ) ]^{1/p} \ \ \mbox{for} \ p \in [1, \infty) .\\[-1mm]
\]
Then $L^{p,\infty} = \{ f : | f |_{p,\infty} < \infty \}$, and it is clear that $L^p
\subset L^{p,\infty}$ for $p \in [1, \infty)$, because $| f |_{p,\infty} \leq \| f
\|_{p}$ in view of Chebyshev's inequality; the Marcinkiewicz space for $p = \infty$
is $L^\infty$ by definition.

It occurs that $| f |_{p,\infty}$ is not a norm for $p \in [1, \infty)$, but a
quasi-norm because\\[-1mm]
\[ | f + g |_{p,\infty} \leq 2 ( | f |_{p,\infty} + | g |_{p,\infty} )\\[-1mm]
\]
(see, e.g., \cite{BL}, p. 7). However, it is possible to endow $L^{p,\infty}$, $p
\in (1, \infty)$, with a norm $\| \cdot \|_{p,\infty}$ converting it into a Banach
space; moreover, the inequality\\[-2mm]
\[ | f |_{p,\infty} \leq \| f \|_{p,\infty} \leq p (p - 1)^{-1} | f |_{p,\infty}\\[-1mm]
\]
holds for all $f \in L^{p,\infty}$. It is worth mentioning that $L^{p,\infty}$
belongs (as a limiting case) to the class of {\it Lorentz spaces} $L^{p,q}$, $q \in
[1, \infty]$ (see, e.g., \cite{BL}, sect.~1.6, and references cited in this book).

Another generalisation of $L^{p,\infty}$, known as the Marcinkie\-wicz space
$M_\varphi$, is defined with the help of a non-negative, concave function $\varphi
\in C [0, \infty)$. This Banach space consists of all (equivalence classes of)
measurable functions for which the following norm\\[-1mm]
\[ \| f \|_\varphi = \sup_{t>0} \frac{1}{\varphi (t)} \int_0^t f^* (s) \, \D s\\[-1mm]
\]
is finite. Here $f^*$ denotes the non-increasing rearrangement of $f$, i.e.,
\[ f^* (s) = \inf_{\lambda > 0} \{ \lambda : m ( \{ x : | f (x) | > \lambda \} ) 
\leq s \} \ \ \mbox{for} \ s \geq 0 ,
\]
and so is non-negative and right-continuous; moreover, its distribution function $m
( \{ x : | f^* (x) | > \lambda \} )$ coincides with that of $f$. If $\varphi (t) =
t^{1 - 1/p}$, then the corresponding Marcin\-kiewicz space is $L^{p,\infty}$,
whereas $\varphi (t) \equiv 1$ and $\varphi (t) = t$ give $L^1$ and $L^\infty$,
respectively.

{\bf The Marcinkiewicz interpolation theorem for boun\-ded linear operators.} This
kind of continuous operators is usually considered as mapping one normed space to
another one, in which case the operator's norm is an important characteristic.
However, the latter can be readily generalised for a mapping of $L^{p}$ to
$L^{p,\infty}$; indeed, if $| T f |_{p,\infty} \leq C \| f \|_{p}$, then it is
natural to introduce the norm (or quasi-norm) of $T$ as the infimum over all
possible values of $C$. Now we are in a position to formulate.

\vspace{1mm}

\noindent {\bf Theorem 2.} {\it Let $p_0, \ q_0, \ p_1, \ q_1 \in [1, \infty]$
satisfy the inequalities $p_0 \leq q_0$, $p_1 \leq q_1$ and $q_0 \neq q_1$, and let
$p, \ q \in [1, \infty]$ be such that $p \leq q$ and the equalities
\[ \frac{1}{p} = \frac{1 - \theta}{p_0} + \frac{1}{p_1} \ \ and \ \ \frac{1}{q} = 
\frac{1 - \theta}{q_0} + \frac{1}{q_1}
\]
hold for some $\theta \in (0, 1)$. If $T$ is a linear operator, which maps $L^{p_0}$
into $L^{q_0,\infty}$ and its norm is $N_0$ and simultaneously $T : L^{p_1} \to
L^{q_1,\infty}$ has $N_1$ as its norm, then $T$ maps $L^{p}$ into $L^{q}$ and its
norm $N$ satisfies the estimate\\[-1mm]
\begin{equation}
N \leq C N_0^{1 - \theta} N_1^\theta\\[-1mm] \label{conv}
\end{equation}
with $C$ depending on $p_0, \ q_0, \ p_1, \ q_1$ and $\theta$.}

\vspace{1mm}

The convexity inequality \eqref{conv} is a characteristic feature of the
interpolation theory. The general form of this theorem (it is valid for
quasi-additive operators, whose special case are sub-additive ones described prior
to Theorem~1) is proved in \cite{ZS}, ch.~XII, sect.~4. In particular, it is shown
that one can take
\[ C = 2 \left( \frac{q}{|q-q_0|} + \frac{q}{|q-q_1|} \right)^{1/q} \frac{p_0^{(1-\theta)
/p_0} p_1^{\theta /p_1}}{p^{1/p}} \, ;
\]
see \cite{ZS}, vol.~II, p.~114, formula (4.18), where, unfortunately, the notation
differs from that adopted here. Special cases of Theorem~2 and diagrams illustrating
them can be found in \cite{LM1}, pp.~155--156. It should be emphasised that the
restriction $p \leq q$ is essential; indeed, as early as 1964, R.\,A.~Hunt \cite{Hu}
constructed an example demonstrating that Theorem 2 is not true without it; for a
description of this example see, e.g., \cite{BL}, pp.~16--17.

It was Marcinkiewicz himself who proposed an extension of his interpolation theorem
to other function spaces. Namely, the so-called diagonal case (when $p_0 = q_0$ and
$p_1 = q_1$) of his theorem is formulated for Orlicz spaces in \cite{M3}. References
to papers containing further results on interpolation in these and other spaces
(e.g., Lorentz and $M_\varphi$) can be found in \cite{BL}, pp.~128--129, and
\cite{LM1}, pp.~163--166.

{\bf Applications of the interpolation theorems.} (1) In his monograph \cite{ZS},
Zygmund gave a detailed study of the one-dimensional Fourier transform:\\[-1mm]
\[ F (f) (\xi) = \frac{1}{\sqrt{2 \pi}} \int_{\RR} f (x) \, \exp 
\{-i \, \xi x\} \, \D x , \ \ \xi \in \RR .
\]
See vol.~II, ch.~XVI, sections~2 and 3, where, in particular, it is demonstrated
that $F$, originally defined on a dense set in $L^p$, $p \in [1, 2]$, is extensible
to the whole space as a bounded operator $F : L^p \to L^{p'}$, $p' = p / (p-1)$, and
so the integral converges in $L^{p'}$. To prove this assertion and its
$n$-dimensional analogue one can use Theorem~2; indeed, $F: L^1 \to L^\infty$ is
bounded (this is straightforward to see), and by Plancherel's theorem $F$ is bounded
on $L^2$, and so this theorem is applicable. On the other hand, the Riesz--Thorin
theorem, which has no restriction $p \leq q$, yields a more complete result valid
for the inverse transform $F^{-1}$ as well. The latter operator acting from $L^{p'}$
to $L^p$ is bounded; here $p' \in [2, \infty)$, and so $p = p' / (p'-1) \in (1,
2]$.

(2) In studies of conjugate Fourier series, the following singular integral operator
(the periodic Hilbert transform)\\[-1mm]
\[ H (f) (s) = \frac{1}{2 \pi} \lim_{\epsilon \to 0} \int_{\epsilon \leq |t| \leq
\pi} f (s - t) \cot \frac{t}{2} \, \D t
\]
plays an important role. Indeed, by linearity it is sufficient to define $H$ on a
basis in $L^2 (-\pi, \pi)$, and the relations\\[-1mm]
\[ H (\cos n t) = \sin n s \ \mbox{for} \ n \geq 0 , \ \ H (\sin n t)
= - \cos n s \ \ \mbox{for} \ n \geq 1\\[-1mm]
\]
show that it expresses passing from a trigonometric series to its conjugate.
Moreover, these formulae show that $H$ is bounded on $L^2 (-\pi, \pi)$ and its norm
is equal to one.

In the mid-1920s, Marcel Riesz obtained his celebrated result about this operator;
first, he announced it in a brief note in the {\it Comptes rendus de l'Acad\'emie
des Sciences Paris}, and three years later published his rather long proof that $H$
is bounded on $L^p (-\pi, \pi)$ for $p \in (1, \infty)$, i.e. for every finite $p >
1$ there exists $A_p > 0$ such that\\[-1mm]
\begin{equation}
\| H (f) \|_p \leq A_p \| f \|_p \ \ \mbox{for all} \ f \in L_p (-\pi, \pi)
.\\[-1mm]
\label{MR}
\end{equation}
However, \eqref{MR} does not hold for $p = 1$ and $\infty$; see \cite{ZS}, vol.~I,
ch.~VII, sect.~2, for the corresponding examples and a proof of this inequality.

There are several different proofs of this theorem; the original proof of M. Riesz
was reproduced in the first edition of the Zygmund's monograph \cite{ZS} which
appeared in 1935. In the second edition published in 1959, this proof was replaced
by that of Calder\'on obtained in 1950. Let us outline another proof based on the
Marcinkiewicz interpolation theorem analogous to Theorem~1, but involving
$L^p$-spaces on $(-\pi, \pi)$ instead of the spaces on $\RR$.

First we notice that it is sufficient to prove \eqref{MR} only for $p \in (1, 2]$.
Indeed, assuming that this is established, then for $f \in L_p$ and $g \in L_{p'}$
we have\\[-2mm]
\[ \int_{-\pi}^\pi [H (f) (s)] \, g (-s) \, \D s \leq A_p \| f \|_p \, \| g 
\|_{p'} \, .\\[-1mm]
\]
by the H\"older inequality (as above $p' = p / (p-1)$, and so $p' \geq 2$ when $p
\leq 2$). Since\\[-2mm]
\[ \int_{-\pi}^\pi [H (f) (s)] \, g (-s) \, \D s = \int_{-\pi}^\pi f (-s) \, [H (g) 
(s)] \, \D s  \, ,\\[-1mm]
\]
the inequality $\| H (g) \|_{p'} \leq A_p^{-1} \| g \|_{p'}$ is a consequence of the
assertion converse to the H\"older inequality.

It was mentioned above that $H$ is bounded in $L^2$. Hence, in order to apply
Theorem~1 for $p \in (1, 2]$, it is sufficient to show that this operator is of weak
type $(1, 1)$, and this is an essential part of Calder\'on's proof; see \cite{ZS},
vol.~I, ch.~IV, sect.~3. Moreover, an improvement of the latter proof allowed
S.\,K.~Pichorides \cite{P} to obtain the least value of the constant $A_p$ in
\eqref{MR}. It occurs that $A_p = \tan \pi / (2 p)$ and $\cot \pi / (2 p)$ is this
value for $p \in (1, 2]$ and $p \geq 2$, respectively.

There are many other applications of interpolation theorems in analysis; see, e.g.,
\cite{BL}, ch.~1, \cite{ZS}, ch.~XII, and references cited in these books.

{\bf Further development of interpolation theorems.} Results constituting the
interpolation space theory were obtained in the early 1960's and are classical now.
This theory was created in the works of Nachman Aronszajn, Alberto Calder\'on,
Mischa Cotlar, Emilio Gagl\-iardo, Selim Grigorievich Krein, Jacques-Louis Lions and
Jaak Peetre to list a few. We leave aside several versions of complex interpolation
spaces developed from the Riesz--Thorin theorem (see, e.g., \cite{BL}, ch.~4), and
concentrate on `espaces de moyennes' introduced by Lions and Peetre in their
celebrated article \cite{LP}. These `real method interpolation spaces' usually
denoted $(A_0, A_1)_{\theta,p}$ are often considered as `descendants' of the
Marcinkiewicz interpolation theorem.

Prior to describing these spaces, it is worth mentioning another germ of
interpolation theory originating from Lw\'ow. The problem~87 in the {\it Scottish
Book} \cite{SB} posed by Banach demonstrates his interest in nonlinear
interpolation. Presumably, it was formulated during Marcinkiewicz's stay in Lw\'ow;
indeed, he solved problems 83 and 106 in \cite{SB}, which were posed before and
after, respectively, the Banach's problem on interpolation. A positive solution of
the latter problem (due to L.~Maligranda) is presented in \cite{SB}, pp.~163--170.

Let us turn to defining the family of spaces $\{ (A_0, A_1)_{\theta,p} \}$ involved
in the real interpolation method; here $\theta \in (0, 1)$ and $p \in [1, \infty]$.
In what follows, we write $A_{\theta,p}$ instead of $(A_0, A_1)_{\theta,p}$ for the
sake of brevity. Let $A_0$ and $A_1$ be two Banach spaces, both continuously
embedded in some (larger) Hausdorff topological vector space, then for a pair
$(\theta,p)$ the space $A_{\theta,p}$ with $p < \infty$ consists of all $a \in A_0 +
A_1$ for which the following norm\\[-2mm]
\[ \| a \|_{\theta,p} = \left\{ \int_0^\infty \Big[ t^{-\theta} \, K (t, a) \Big]^p
\frac{\D t}{t} \right\}^{1/p}\\[-1mm]
\]
is finite. Here $K (t, a)$ is defined on $A_0 + A_1$ for $t \in (0, \infty)$
by\\[-2mm]
\[ \inf_{a_0, a_1} \{ \| a_0 \|_{A_0} + t \| a_1 \|_{A_1} : a_0 \in A_0 , 
\ a_1 \in A_1 \ \mbox{and} \ a_0 + a_1 = a \} .\\[-2mm]
\]
This $K$-functional was introduced by Peetre. If $p = \infty$, then the expression
$\sup_{t > 0} \{ t^{-\theta} \, K (t, a) \}$ gives the norm $\| a \|_{\theta,
\infty}$ when finite.

Every $A_{\theta,p}$ is an intermediate space with respect to the pair $(A_0, A_1)$,
i.e.,\\[-2mm]
\[ A_0 \cap A_1 \subset A_{\theta,p} \subset A_0 + A_1 .\\[-1mm]
\]
Moreover, if $A_0 \subset A_1$, then\\[-2mm]
\[ A_0 \subset A_{\theta_0,p_0} \subset A_{\theta_1,p_1} \subset A_1\\[-1mm]
\]
provided either $\theta_0 > \theta_1$ or $\theta_0 = \theta_1$ and $p_0 \leq p_1$.
For any $p$, it is convenient to put $A_{0,p} = A_0$ and $A_{1,p} = A_1$. Now we are
in a position to explain what the interpolation of an operator is in terms of the
family $\{ A_{\theta,p} \}$ and another family of spaces $\{ B_{\theta,p} \}$
constructed by using some Banach spaces $B_0$ and $B_1$ in the same way as $A_0$
and $A_1$.

Let $T : A_0 + A_1 \to B_0 + B_1$ be a linear operator such that its norm as the
operator mapping $A_0 \, (A_1)$ to $B_0 \, (B_1)$ is equal to $M_0 \, (M_1)$, then
the operator $T : A_{\theta,p} \to B_{\theta,p}$ is also bounded and its norm is
less than or equal to $M_0^{1 - \theta} M_1^\theta$. Along with the method based on
the $K$-functional, there is an equivalent method (also developed by Peetre)
involving the so-called $J$-functional. Further details concerning this approach to
interpolation theory can be found in \cite{BL}, chapters~3 and 4.

\subsection*{The Marcinkiewicz function}

In the {\it Annales de la Soci\'et\'e Polonaise de Math\'ematique}, volume 17
(1938), Marcinkiewicz published two short papers. Two remarkable integral operators
were considered in the first of these notes (see \cite{M4} and \cite{JM},
pp.~444--451); they and their numerous generalisations became indispensable tools in
analysis. One of these operators is always called the `Mar\-cinkiewicz integral';
see \cite{ZS}, ch.~IV, sect.~2, for its definition and properties. In particular, it
is used for investigation of the structure of a measurable set near `almost
arbitrary' point; see \cite{S}, sections 2.3 and 2.4, whereas further references to
papers describing some its generalisations can be found in the monographs \cite{S}
and \cite{ZS}. The second operator is usually referred to as the `Marcinkiewicz
function' (see, e.g., \cite{LM1}, pp.\,192--194), but it also appears as the
`Marcinkiewicz integral'. Presumably, the mess with names began as early as 1944,
when Zygmund published the extensive article \cite{Z1}, section 2 of which was
entitled ``On an integral of Marcinkiewicz''. In fact, this 14-pages long section is
devoted to a detailed study of the Marcinkiewicz function $\mu$, whose properties
were just outlined by Marcinkiewicz himself in \cite{M4}. It is not clear whe\-ther
Zygmund had already received information about Mar\-cinkiewicz's death, when he
decided to present in detail the results from \cite{M4} (the discovery of mass
graves in the Katy\'n forest was announced by the Nazi government in April 1943).

Zygmund begins his presentation with a definition of the Littlewood--Paley function
$g (\theta; f)$, which is a nonlinear operator applied to an integrable, $2
\pi$-periodic $f$. The purpose of introducing $g (\theta; f)$ was to provide a
characterisation of the $L^p$-norm $\| f \|_p$ in terms of the Poisson integral of
$f$. After describing some properties of $g (\theta)$, Zygmund notes.\\[-4mm]
\begin{quote}
It is natural to look for functions analogous to $g (\theta)$ but defined without
entering the interior of the unit circle.\\[-4mm]
\end{quote}
After a reference to \cite{M4}, Zygmund continues:\\[-4mm]
\begin{quote}
Marcinkiewicz had the right idea of introducing the function\\[-8mm]
\end{quote}
\begin{eqnarray*}
&& \!\!\!\!\!\!\!\!\!\! \mu (\theta) = \mu (\theta; f) \\ && = \Big\{ \int_0^{\pi}
\frac{[F (\theta + t) + F (\theta - t) - 2 F (\theta)]^2}{t^3} \D t \Big\}^{1/2} \\
&& = \Big\{ \int_0^{\pi} t \Big[ \frac{F (\theta + t) + F (\theta - t) - 2 F
(\theta)}{t^2} \Big]^2 \D t \Big\}^{1/2}
\end{eqnarray*}
\begin{quote}
where $F (\theta)$ is the integral of $f$,\\[-8mm]
\end{quote}
\[ F (\theta) = C + \int_0^\theta f (u) \, \D u \, .
\]
\begin{quote}
More generally, he considers the functions\\[-8mm]
\end{quote}
\begin{eqnarray*}
\ \ \ \ \, \mu_r (\theta) = \Big\{ \int_0^{\pi} \frac{|F (\theta + t) + F (\theta -
t) - 2 F (\theta)|^r}{t^{r+1}} \D t \Big\}^{1/r} \\ = \Big\{ \int_0^{\pi} t^{r-1}
\left| \frac{F (\theta + t) + F (\theta - t) - 2 F (\theta)}{t^2} \right|^r \D t
\Big\}^{1/r} ,
\end{eqnarray*}
\begin{quote}
so that $\mu_2 (\theta) = \mu (\theta)$. He proves the following facts which are
clearly analogues of the corresponding properties of $g (\theta)$.\\[-6mm]
\end{quote}
These facts are the estimates\\[-1mm]
\[ \| \mu_q \|_q \leq A_q \| f \|_q \ \ \mbox{and} \ \ \| f \|_p \leq A_p
\| \mu_p \|_p%\\[-1mm]
\]
valid for $q \geq 2$ and $1 < p \leq 2$, respectively, where $f$ has the zero mean
value in the second inequality, and the assertion: {\it For every $p \in (1, 2]$
there exists a continuous, $2 \pi$-periodic function $f$ such that $\mu_p (\theta; f)
= \infty$ for almost every $\theta$.}

Furthermore, Marcinkiewicz conjectured that for $p > 1$ the inequalities\\[-2mm]
\begin{equation}
A_p \| f \|_p \leq \| \mu \|_p \leq B_p \| f \|_p\\[-1mm] \label{m-z}
\end{equation}
hold, where again $f$ must have the zero mean value in the second inequality.
Moreover, he foresaw that it would not be easy to prove these inequalities; indeed,
the proof given by Zygmund in his article \cite{Z1} is more than 11 pages long.

The first step towards generalisation of the Marcinkie\-wicz function was made by
Daniel Waterman; his paper \cite{W} was published seven (!) years after presentation
of the work to the AMS. However, its abstract appeared in the {\it Proceedings of
the International Congress of Mathema\-ticians} held in 1954 in Amsterdam. Waterman
considered the following $\mu$-function\\[-1mm]
\[ \mu (\tau; f) = \Big\{ \int_0^{\infty} \frac{[F (\tau + t) + F (\tau - t) - 2 F
(\tau)]^2}{t^3} \D t \Big\}^{1/2} ,%\\[-1mm]
\]
where $\tau \in (-\infty, \infty)$ and $F$ is a primitive of $f \in L^p (-\infty,
\infty)$, $p>1$. His proof of inequalities \eqref{m-z} for $\mu (\tau; f)$ heavily
relies on the M.~Riesz theorem about conjugate functions on $\RR^1$ (see \cite{W},
p.~130, for the formulation), and its proof involves the Marcinkiewicz
interpolation theorem described above.

Another consequence of inequalities \eqref{m-z} for $\mu (\tau; f)$ is a
characterization of the Sobolev space $W^{1,p} (\RR)$, $p \in (1, \infty)$. Indeed,
putting\\[-1mm]
\[ \mathrm{M} (\tau; f) = \Big\{ \int_0^{\infty} \frac{[f (\tau + t) + f (\tau - t) 
- 2 f (\tau)]^2}{t^3} \D t \Big\}^{1/2}
\]
for $f \in W^{1,p} (\RR)$, we have that $\mathrm{M} (\tau; f) = \mu (\tau; f')$.
Then \eqref{m-z} can be written as\\[-2mm]
\[ A_p \| f' \|_p \leq \| \mathrm{M} (\cdot; f) \|_p \leq B_p \| f' \|_p \, ,\\[-1mm] \]
which implies the following assertion. {\it Let $p \in (1, \infty)$, then $f \in
W^{1,p} (\RR)$ if and only if $f \in L^{p} (\RR)$ and\/ $\mathrm{M} (\cdot; f) \in
L^{p} (\RR)$.}

Stein extended these results to higher dimensions in the late 1950s and early 1960s
(it is worth mentioning that $\mu$ is referred to as the Marcinkiewicz integral in
his paper \cite{S1}). For this purpose he applied the real-variable technique used
in the generalisation of the Hilbert transform\\[-1mm]
\[ {\rm P.V.} \int_0^{\infty} \frac{f (x + t) - f (x - t)}{t} \, \D t \, ,\\[-1mm]
\]
to higher dimensions. Indeed, this can be written as\\[-1mm]
\[ \int_0^{\infty} \frac{F (x + t) + F (x - t) - 2 F (x)}{t^2} \, \D t\\[-1mm]
\]
which resembles the expression for $\mu (\tau; f)$, and so Stein, in his own words,
was\\[-5mm]
\begin{quote}
guided by the techniques used by A.\,P. Calder\'on and A. Zygmund \cite{CZ} in their
study of the $n$- di\-men\-sional generalizations of the Hilbert trans\-form;
connected with this are some earlier ideas of Marcinkiewicz.\\[-5mm]
\end{quote}
The definition of singular integral given in \cite{CZ}, to which Stein refers,
involves a function $\Omega (x)$ defined for $x \in \RR^n$ and assumed: (i) to be
homogeneous of degree zero, i.e. to depend only on $x' = x / |x|$; (ii) to satisfy
the H\"older condition with exponent $\alpha \in (0, 1]$; (iii) to have the zero
mean value over the unit sphere in $\RR^n$. Then\\[-1mm]
\[ S (f) (x) = \lim_{\epsilon \to 0} \int_{|y|>\epsilon} \frac{\Omega (y')}{|y|^n} \, 
f (x-y) \, \D y\\[-1mm]
\]
exists almost everywhere provided $f \in L^p (\RR^n)$, $p \in [1, \infty)$.
Furthermore, this singular integral operator is bounded in $L^p (\RR^n)$ for $p >
1$, i.e. the inequality $\| S (f) \|_p \leq A_p \| f \|_p$ holds with $A_p$
independent of $f$.

Moreover, in the section dealing with background facts, Stein notes that $\mu$ is a
nonlinear operator and writes (see \cite{S1}, p.~433):\\[-5.6mm]
\begin{quote}
An ``interpolation'' theorem of Marcinkiewicz is very useful in this connection.

In quoting the result of Marcinkiewicz, [\dots]~we shall not aim at generality. For
the sake of~sim\-plicity we shall limit ourselves to the special case that is
needed.\\[-5mm]
\end{quote}
After that the required form of the interpolation theorem (see Theorem~1 above) is
formulated and used later in the paper, thus adding one of the first items in the
now long list of its applications. Since the term interpolation was novel, quotation
marks are used by Stein in the quoted piece; indeed, Zygmund's proof of the
Marcinkiewicz theorem had appeared in 1956, just two years earlier than Stein's
article.

His generalization of the Marcinkiewicz function $\mu (\tau; f)$ Stein begins with
the case when $f \in L^p (\RR^n)$, $p \in [1, 2]$. Realising the analogy described
above, he puts
\begin{equation}
F_t (x) = \int_{|y| \leq t} \frac{\Omega (y')}{|y|^{n-1}} \, f (x-y) \, \D y \, , \
\ x \in \RR^n , \label{S1}
\end{equation}
where $\Omega$ satisfies conditions (i)--(iii), and notes that if $n=1$ and $\Omega
(y) = {\rm sign}\, y$, then\\[-2mm]
\[ F_t (x) = F (x + t) + F (x - t) - 2 F (x) \ \mbox{with} \ F (x) = \int_0^x f (s) 
\, \D s .
\]
Therefore, it is natural to define the $n$-dimensional Mar\-cin\-kiewicz function as
follows:\\[-2mm]
\begin{equation}
\mu (x; f) = \Big\{ \int_0^{\infty} \frac{[F_t (x)]^2}{t^3} \, \D t \Big\}^{1/2}
.\\[-1mm] \label{S2}
\end{equation}

His investigation of properties of this function Stein begins by proving that $\|
\mu (\cdot; f) \|_2 \leq A \| f \|_2$, where $A$ is independent of $f$, and his
proof involving Plancherel's theorem is not elementary at all. Even less elementary
is his proof that $\mu (\cdot; f)$ is of weak type $(1, 1)$. Then the Marcinkiewicz
interpolation theorem (see Theorem~1 above) implies that $\| \mu (\cdot; f) \|_p
\leq A \| f \|_p$ for $p \in (1, 2]$ provided $f \in L^p (\RR^n)$. For all $p \in
(1, \infty)$ this inequality is proved in \cite{S1} with assumptions (i)--(iii)
changed to the following ones: $\Omega (x')$ is absolutely integrable on the unit
sphere and is odd there, i.e. $\Omega (-x') = - \Omega (x')$. A few years later,
A.~Benedek, A.\,P. Cal\-der\'on and R.~Panzone demonstrated that for a
$C^1$-function $\Omega$ condition (iii) implies the last inequality for all $p \in
(1, \infty)$.

In another note, Stein obtained the following generalisation of the one-dimensional
result. {\it Let $p \in (2 n / (n+2), \infty)$ and $n \geq 2$, then $f$ belongs to
the Sobolev space $W^{1,p} (\RR^n)$ if and only if $f \in L^{p} (\RR^n)$
and}\\[-2mm]
\[ \Big\{ \int_{\RR^n} \frac{[f (\cdot + y) + f (\cdot - y) - 2 f 
(\cdot)]^2}{|y|^{n+2}} \D y \Big\}^{1/2} \in L^{p} (\RR^n) .\\[-1mm]
\]
For $n > 2$ this does not cover $p \in (1, 2 n / (n+2)]$, and so is weaker than the
assertion formulated above for $n = 1$.

In the survey article \cite{LM1}, pp.~193--194, one finds a list of papers
concerning the Marcinkiewicz function. In particular, further properties of $\mu$
were considered by A.~Torchinsky and S.~Wang \cite{TW} in 1990, whereas T.~Walsh
\cite{Wa} proposed a modification of the definition \eqref{S1}, \eqref{S2} in 1972.

\subsection*{Multipliers of Fourier series and integrals}

During his stay in Lw\'ow, Marcinkiewicz collaborated with Stefan Kaczmarz and
Juliusz Schauder,$^*$ \footnote{$^*$Both perished in World War II. Being in the
reserve, Kaczmarz was drafted and killed during the first week of war; the
circumstances of his death are unclear. Schauder was on hiding in occupied Lw\'ow
and the Gestapo killed him in 1943 while he was trying to escape arrest.} due to
whom~his interest in {\it multipliers} of orthogonal series had arisen. Studies in
this area of analysis were initiated by Hugo Steinhaus in the 1920s; in its general
form, the problem of multipliers is as follows. Let $B_1$ be a Banach space with a
Schauder basis $\{ g_n \}_{n=1}^\infty$, the (linear) operator $T$ is called
multiplier when there is a sequence $\{ m_n \}_{n=1}^\infty$ of scalars of this
space and $T$ acts as follows:\\[-2mm]
\[ B_1 \ni f = \sum_{n=1}^\infty c_n g_n \ \ \to \ \ T f \sim \sum_{n=1}^\infty m_n
c_n g_n \, .\\[-2mm]
\]
Here $\sim$ means that the second sum assigned as $T f$ can belong to the same space
$B_1$ or be an element of another Banach space $B_2$; this depends on properties of
the sequence. Multipliers of Fourier series are of paramount interest and this was
the topic of the remarkable paper \cite{M1} published by Marcinkiewicz in 1939.

Not long before Marcinkiewicz's visit to Lw\'ow started, Kaczmarz investigated some
properties of multipliers in the function spaces (mainly $L^p (0, 1)$ and $C [0,
1]$) under rather general assumptions about the system $\{ g_n \}_{n=1}^\infty$.
Further results about multiplier operators were obtained in the joint paper
\cite{KM} of Kaczmarz and Marcinkiewicz; it was submitted to the {\it Studia
Mathematica} in June 1937, i.e., their collaboration lasted for another year after
Marcinkiewicz left Lw\'ow. This paper has the same title as that of Kaczmarz and
concerns the case when $L^p (0, 1)$ with $p \neq \infty$ is mapped to $L^q (0, 1)$,
$q \in [1, \infty]$; it occurs that the case $q=\infty$ is the simplest one. In this
paper, it is assumed that every function $g_n$ is bounded, whereas the sequence $\{
g_n \}_{n=1}^\infty$ is closed in $L^1 (0, 1)$. In each of four theorems which
differ by the ranges of $p$ and $q$ involved, certain conditions are imposed on $\{
m_n \}_{n=1}^\infty$ and these conditions are necessary and sufficient for the
sequence to define a multiplier operator $T: L^p \to L^q$.

After returning to Wilno, Marcinkiewicz kept on his studies of multipliers initiated
in Lw\'ow, and in May 1938, he submitted (again to the {\it Studia Mathematica}) the
seminal paper \cite{M1}, in which the main results are presented in a curious way.
Namely, Theorems~1 and 2, concerning multipliers of Fourier series and double
Fourier series, are formulated in the reverse order. Presumably, the reason for this
is the importance of multiple Fourier series for applications and generalisations.
Let us formulate Theorem~1 in a slightly updated form.

\vspace{1mm}

{\it Let $f \in L^p (0, 2 \pi)$, $p \in (1, \infty)$, be a real-valued function and
let its Fourier series be\\[-2mm]
\[ a_0 / 2 + \sum_{n=1}^\infty A_n (x) , \ where \ A_n (x) = a_n \cos n x +
b_n \sin n x .\\[-1mm]
\]
If a bounded sequence $\{ \lambda_n \}_{n=1}^\infty \subset \RR$ is such
that\\[-2mm]
\begin{equation} 
\sum_{n=2^k}^{2^{k+1}} |\lambda_n - \lambda_{n+1}| \leq M \ \ for \ all \ \
k=0,1,2,\dots,\\[-1mm] \label{m-m}
\end{equation}
where $M$ is a constant independent of $k$, then the mapping $f \mapsto
\sum_{n=1}^\infty \lambda_n \, A_n$ is a bounded operator in $L^p (0, 2 \pi)$.}

\vspace{1mm}

\noindent It is well-known that for $p=2$ this theorem is true with condition
\eqref{m-m} omitted, but this is not mentioned in \cite{M1}. The assumptions that
$f$ is real-valued and $\{ \lambda_n \}_{n=1}^\infty \subset \RR$ were not stated in
\cite{M1} explicitly, but used in the proof. This was noted by Solomon Grigorievich
Mikhlin \cite{SM}, who extended this theorem to complex-valued multipliers and
functions; also, he used the exponential from of the Fourier expansion:\\[-3mm]
\[ f (x) = \sum_{n=-\infty}^\infty \!\!\!\! c_{n} \exp i n x \, .\\[-2mm]
\]

The trigonometric form was used by Marcinkiewicz for double Fourier series as well,
and his sufficient conditions on bounded real multipliers $\{ \lambda_{m n} \}$ look
rather awkward. Now, the restrictions on $\{ \lambda_{m n} \} \subset \CC$ are
usually expressed in a rather condensed form by using the so-called dyadic
intervals; see, e.g., \cite{S}, sect.~5.1. Applying these conditions to multipliers
acting on the expansion\\[-2mm]
\[ \sum_{m,n=-\infty}^\infty \!\!\!\! c_{mn} \exp i \{ m x + n y \}\\[-1mm]
\]
of $f \in L^p ( (0, 2 \pi)^2 )$, $p \in (1, \infty)$, one obtains an updated
formulation of the multiplier theorem; see, e.g., \cite{LM1}, p.~201.

A simple corollary derived by Marcinkiewicz from this theorem is as follows (see
\cite{M1}, p.~86). The fractions\\[-2mm]
\begin{equation}
\frac{m^2}{m^2 + n^2} , \ \ \frac{n^2}{m^2 + n^2} , \ \ \frac{|m n|}{m^2 + n^2}
\label{m-n}
\end{equation}
provide examples of multipliers in $L^p$ for double Fourier series. The reason to
include these examples was to answer a question posed by Schauder and this is
specially mentioned in a footnote. Moreover, after remarking that his Theorem~2
admits an extension to multiple Fourier series, Marcinkiewicz added a
straightforward generalisation of formulae \eqref{m-n} to higher dimensions again
referring to Schauder's question. This is an evidence that the question was an
important stimulus for Marcinkiewicz in his work.

A natural way to generalise Marcinkiewicz's theorems is to consider multipliers of
Fourier integrals. Study of these operators was initiated by Mikhlin in 1956; see
note \cite{SM} in which the first result of that kind was announced. Several years
later, Mikhlin's theorem was improved by Lars H\"ormander \cite{H}, and since than
it is widely used for various purposes. To formulate this theorem we need the
$n$-dimen\-sional Fourier transform
\[ F (f) (\xi) = (2 \pi)^{-n/2} \int_{\RR^n} f (x) \, \exp 
\{-i \, \xi \cdot x\} \, \D x , \ \ \xi \in \RR^n ,
\]
defined for $f \in L^2 (\RR^n) \cap L^p (\RR^n)$, $p \in (1, \infty)$. It is clear
that any bounded measurable function $\Lambda$ on $\RR^n$ defines the
mapping\\[-1mm]
\[ T_\Lambda (f) (x) = F^{-1} [ \Lambda (\xi) F (f) (\xi)] (x) \, , \ \ x \in \RR^n 
,
\]
such that $T_\Lambda (f) \in L^2 (\RR^n)$. If $T_\Lambda (f)$ is also in $L^p
(\RR^n)$ and $T_\Lambda$ is a bounded operator, i.e.,\\[-1mm]
\begin{equation}
\| T_\Lambda (f) \|_p \leq B_{p,n} \| f \|_p \ \ \mbox{for all} \ f \in L^p (\RR^n)
\label{mult}
\end{equation}
with $B_p$ independent of $f$, then $\Lambda$ is called a multiplier for~$L^p$.

The description of all multipliers for~$L^2$ is known as well as for $L^1$ and
$L^\infty$ (it is the same for these two spaces); see \cite{S}, pp.~94--95.
However, the question about characterisation of the whole class of multipliers for
other values of $p$ is far from its final solution. The following assertion gives
widely used sufficient conditions.

\vspace{1mm}

\noindent {\bf Theorem} (Mikhlin, H\"ormander). {\it Let $\Lambda$ be a function of
the $C^k$-class in the complement of the origin of $\RR^n;$ here $k$ is the least
integer greater than $n/2$. If there exists $B > 0$ such that\\[-1mm]
\[ |\xi|^\ell \left| \frac{\partial^\ell \Lambda (\xi)}{\partial \xi_{j_1} \partial
\xi_{j_2} \dots \partial \xi_{j_\ell}} \right| \leq B \, , \ \ 1 \leq j_1 < j_2 <
\dots < j_\ell \leq n \, ,
\]
for all $\xi \in \RR^n$, $\ell = 0,\dots,k$ and all possible $\ell$-tuples, then
inequality \eqref{mult} holds, i.e., $\Lambda$ is a multiplier for~$L^p$.}

\vspace{1mm}

In various versions of this theorem, different assumptions are imposed on the
differentiability of $\Lambda$. In particular, H\"ormander \cite{H}, pp.~120--121,
replaced the pointwise inequality for weighted derivatives of $\Lambda$ by a weaker
one involving certain integrals (see also \cite{S}, p.~96). Recently, Loukas
Grafakos and Lenka Slav\'ikov\'a \cite{GS} obtained new sufficient conditions for
$\Lambda$ in the multiplier theorem, thus improving H\"ormander's result. Their
conditions are optimal in a certain sense explicitly described in \cite{GS}.

\vspace{1mm}

\noindent {\bf Corollary.} {\it Every function, which is smooth everywhere except at
the origin and is homogeneous of degree zero, is a Fourier multiplier for~$L^p$.}

\vspace{1mm}

Its immediate consequence is the Schauder estimate\\[-2mm]
\[ \left\| \frac{\partial^2 u}{\partial x_{j_1} \partial
x_{j_2}} \right\|_p \leq C_{p,n} \| \Delta \, u \|_p \, , \ \ 1 \leq j_1 , j_2 \leq
n \, ,\\[-1mm]
\]
valid for $u$ belonging to the Schwartz space of rapidly decaying infinitely
differentiable functions. For this purpose one has to use the equality\\[-1mm]
\[ F \left( \frac{\partial^2 u}{\partial x_{j_1} \partial
x_{j_2}} \right) (\xi) = \frac{\xi_{j_1} \xi_{j_2}}{|\xi|^2} \, F ( \Delta \, u )
(\xi) \, , \ \ 1 \leq j_1 , j_2 \leq n \, ,
\]
and the fact that the function $\xi_{j_1} \xi_{j_2} / |\xi|^2$ is homogeneous of
degree zero.

\vspace{2mm}

\noindent {\bf Acknowledgements.} The author thanks Irina Egorova for the photo of
Marcinkiewicz's plaque and Alex Eremen\-ko, whose comments helped to improve the
original manu\-script.

%\vspace{2mm}

$------------------------$

\vspace{2mm}

\noindent {\bf N.\,Kuznetsov} is a Principal Research Scientist in the Laboratory
for Mathematical Modelling of Wave Phenomena at the Institute for Problems in
Mechanical Engineering, Russian Academy of Sciences, St.\,Petersburg. This
laboratory was founded by him in 1997, and he headed it until 2016.

\vspace{1mm}

\noindent e-mail: nikolay.g.kuznetsov@gmail.com

\end{document}